\documentclass[12pt,twoside,a4paper,reqno]{amsart}
\usepackage{amsmath,amsthm,amssymb,amscd,amsfonts}
\usepackage{enumerate}
\makeatletter
\@namedef{subjclassname@2010}{%
  \textup{2010} Mathematics Subject Classification}
\makeatother
\newtheorem{theorem}{Theorem}[section]
\newtheorem{lemma}[theorem]{Lemma}

\theoremstyle{definition}
\newtheorem{definition}[theorem]{Definition}

\newtheorem{remark}[theorem]{Remark}

\theoremstyle{definition}

\theoremstyle{definition}

\numberwithin{equation}{section}

\title{Pointwise amenability for dual
Banach algebras  }

\author{Mannane Shakeri}
\address{Department of Mathematics, Central Tehran Branch, Islamic Azad University,Tehran, Iran, e-mail: {\tt m.s.chenari@gmail.com}}

\author{Amin Mahmoodi}
\address{Department of Mathematics, Central Tehran Branch, Islamic Azad University, Tehran, Iran, e-mail: {\tt a\_mahmoodi@iauctb.ac.ir}}

\begin{document}
\pagestyle{headings}

\begin{abstract}
We shall develop two notions of pointwise amenability, namely
pointwise Connes amenability and pointwise $w^*$-approximate Connes
amenability, for dual Banach algebras which take the $w^*$-topology
into account. We shall study these concepts for the Banach sequence
algebras $\ell^1(\omega)$ and the weighted semigroup algebras $
\ell^{1}(\mathbb{N}_{\wedge},\omega)$. For a weight $\omega$ on a
discrete semigroup $S$, we shall investigate pointwise
amenability/Connes amenability of $\ell^1(S,\omega)$ in terms of
diagonals.
\end{abstract}

\maketitle Keywords: pointwise amenability, pointwise Connes
amenability, Beurling algebras.

\section{Introduction }
The key concept of amenability for Banach algebras introduced by B.
E. Johnson \cite {13}. The pointwise variant of amenability
introduced by H. G. Dales,  F. Ghahramani and R. J. Loy, however
this appeared formally
 in \cite{2}.

Let $\mathcal{A}$ be a Banach algebra. Then the  projective tensor
product $\mathcal{A}\hat{\otimes}\mathcal{A}$ naturally  is a Banach
$\mathcal{A}$-bimodule and the map
$\pi:\mathcal{A}\hat{\otimes}\mathcal{A}\longrightarrow \mathcal{A}$
defined  by  $\pi(a\otimes b)=ab$,  $a,b \in \mathcal{A}$ is a
linear continuous
 $\mathcal{A}$-bimodule homomorphism.
Let $E$ be a Banach $\mathcal{A}$-bimodule. A \textit{derivation} is
a bounded linear map
 $D: \mathcal{A}\longrightarrow E$
 satisfying
$D(ab)=Da\cdot b+a\cdot Db\;\; (a,b\in \mathcal{A})$.
 A Banach algebra  $\mathcal{A}$  is \textit{pointwise  amenable} at $a_0\in
\mathcal{A}$ if, for each  Banach $\mathcal{A}$-bimodule $E$, every
derivation $D: \mathcal{A}\to E^{*}$ is\textit{ pointwise inner } at
$a_0$, that is, there exist $\phi\in E^{*}$ such that
$D(a_0)=a_0\cdot \phi-\phi\cdot a_0$ \cite{2}.

Let $\mathcal{A}$ be a Banach algebra. A Banach
$\mathcal{A}$-bimodule $E$ is \textit{dual} if there is a closed
submodule $E_*$ of $E^*$ such that $E=(E_*)^*$.
 We call $E_*$ the \textit{predual} of $E$. A Banach algebra $\mathcal{A}$ is \textit{dual} if it is dual as a Banach $\mathcal{A}$-bimodule.
 We write $\mathcal{A}=(\mathcal{A}_*)^*$
 if we wish to stress that $\mathcal{A}$ is a dual Banach algebra with predual
 $\mathcal{A}_*$. For a dual Banach algebra $\mathcal{A}$, a dual Banach $\mathcal{A}$-bimodule $E$ is \textit{normal} if the module actions of
  $\mathcal{A}$ on $E$ are
 $w^*$-continuous. The notion of Connes amenability for dual Banach algebras, which is
another modification of the notion amenability
  systematically introduced   by  V. Runde  \cite{9}.
 A dual Banach algebra $\mathcal{A}$ is \textit{Connes
amenable} if every $w^*$-continuous derivation from $\mathcal{A}$
into a normal, dual Banach
 $\mathcal{A}$-bimodule is inner.

  The concept of $w^*$-approximately Connes amenability introduced by the first
 author in \cite{7}. One may see also \cite{3,5,6,8,10,11}, for more information on Connes
 amenability and other related notions.

 The purpose of this paper is to study pointwise Connes amenability
 of dual Banach algebras as well as their pointwise $w^{*}$-approximate Connes amenability. The organization of the paper is as
 follows. In section 2, some basic properties are given. It is
shown that every commutative pointwise Connes amenable dual Banach
algebras must be unital.

In section 3, it is proved that the Banach sequence algebra
$\ell^1(\omega)$ is pointwise $w^*$-approximately Connes amenable
while it is not pointwise Connes amenable, where $\omega$ is a
weight function. It is also shown that the same is true for the
class of weighted semigroup algebras of the form
$\ell^{1}(\mathbb{N}_{\wedge},\omega)$, provided $\lim_n
\omega(n)=\infty$.

In section 4, the relation between pointwise amenability/Connes
amenability of weighted semigroup algebras $\ell^1(S,\omega)$ and
the existence of some specified diagonals is studied. For a discrete
group $G$, the special case $\ell^1(G,\omega)$ is also considered.

In this note, we introduce and study the pointwise variants of
Connes amenability and $w^{*}$-approximate Connes  amenability. We
show that every commutative pointwise Connes amenable dual Banach
algebras must be unital (Theorem \ref{2.5}). For a weight $\omega$
on $\mathbb{N}$, we describe the Banach sequence algebra
$\ell^1(\omega)$. We show that $\ell^1(\omega)$ is pointwise
$w^*$-approximately Connes amenable, while it is not pointwise
Connes amenable (Theorem \ref{3.1}). Indeed, $\ell^1(\omega)$
clearly shows the difference between pointwise $w^*$-approximate
Connes amenability and both pointwise Connes amenability and
$w^*$-approximate Connes amenability. We consider the class of
weighted semigroup algebras of the form
$\ell^{1}(\mathbb{N}_{\wedge},\omega)$, where $ \omega : \mathbb{N}
\longrightarrow [1,\infty)$ is any function such that $\lim_n
\omega(n)=\infty$. We show that $
\ell^{1}(\mathbb{N}_{\wedge},\omega)$ is not pointwise Connes
amenable, however it is pointwise $w^*$-approximately Connes
amenable (Theorem \ref{3.2}).

\section{Connes amenability; pointwise versions }

From \cite{7}, we recall that a dual Banach algebra $\mathcal{A}$ is
\textit{$w^{*} $-approximately Connes  amenable}   if, for every
normal, dual Banach
 $\mathcal{A}$-bimodule $E$, every  $w^{*}$-continuous derivation $D: \mathcal{A}\to E$ is   \textit{$w^{*} $-approximately inner},
 that is, there exists a net $(x_{n})\subseteq E $ such that
           $D(a)$=$w^{*}$-$\lim_{n}(a\cdot \phi_{n}-\phi_{n}\cdot a)$.

We first introduce the pointwise versions of ($w^*$-approximate)
Connes amenability.
\begin{definition}\label{2.1}
Let $\mathcal{A}$ be a dual Banach algebra. Then:

$(i)$ $\mathcal{A}$ is \textit{pointwise Connes amenable} at $a_0\in
\mathcal{A}$ if for every normal, dual Banach $\mathcal{A}$-bimodule
$E$, every $w^*$-continuous derivation $D: \mathcal{A}\to E$ is
pointwise inner at $a_0$.

 $(ii)$ $\mathcal{A}$ is \textit{pointwise
$w^*$-approximately Connes amenable}  at $a_0\in \mathcal{A}$ if for
every normal, dual Banach $\mathcal{A}$-bimodule $E$, every
$w^*$-continuous derivation $D: \mathcal{A}\to E$ is
\textit{pointwise $w^*$-approximately inner} at $a_0$, that is,
there exists a net $(x_n)\subseteq E$ such that
$D(a_0)=w^*-\lim_n\left(a_0\cdot x_n-x_n\cdot a_0\right)$.

$(iii)$ $\mathcal{A}$ is \textit{ pointwise ($w^*$-approximately)
Connes amenable} if $\mathcal{A}$ is pointwise ($w^*$-approximately)
Connes amenable at each $a\in A$.
\end{definition}

Let $\mathcal{A}$ be a Banach algebra. From \cite{4}, we recall that
$\mathcal{A}$ has \textit{ left (right) approximate units} if, for
each $a\in\mathcal{A}$ and  $\varepsilon>0$ there exists
$u\in\mathcal{A}$ such that  $ \|a-ua\|<\epsilon$  ($
\|a-au\|<\epsilon$),
 and $\mathcal{A}$  has  approximate units if, for each $a\in\mathcal{A}$ and  $\varepsilon>0$, there exists $u\in\mathcal{A}$  such that
  $ \|a-ua\|+\|a-au\|<\epsilon$.
The appropriate approximate units have bound $ m$ if the element  $
u $ can be chosen such that $\|u\|\leq m$. The algebra $\mathcal{A}$
 has a \textit{ bounded} (left or right) approximate units if it has (left or right) approximate units of bound    $ m$    for some
  $ m\geq1$.
\begin{definition}\label{2.2}
A dual Banach algebra  $\mathcal{A}=(\mathcal{A}_*)^*$ has \textit{
left  (right)  $w^{*}$-approximate units} if, for each
$a\in\mathcal{A}$ and  $\varepsilon>0$, and for each finite subset
$\mathcal{K}\subseteq \mathcal{A}_*$,  there is $u\in\mathcal{A}$
such that    $|\langle\psi,a-ua\rangle|<\varepsilon$
($|\langle\psi,au-a\rangle|<\varepsilon$) for
 $\psi\in \mathcal{K}$. We say   $\mathcal{A}$   has    $w^{*}$ -approximate units if, for each    $a\in\mathcal{A}$     and
   $\varepsilon>0$   and for each  finite subset    $\mathcal{K}\subseteq \mathcal{A}_*$, there is     $u\in\mathcal{A}$
    such that  $$|<\psi,a-au>|+|<\psi,a-ua>|<\epsilon \ ,    \ \ (\psi\in \mathcal{K}) \ .$$
The appropriate $w^{*}$-approximate units have bounded $ m$   if the
element $u $ can be chosen such that $ \|u\|\leq m $. The dual
Banach algebra  $\mathcal{A}$ has  \textit{bounded} (left or right)
$w^{*}$-approximate units if it  has (left or right)
 $w^{*}$-approximate units of bound $ m$  for some $ m\geq 1$.
\end{definition}
The following lemma is useful in considerations of identities.
\begin{lemma}\label{2.3}
Let   $\mathcal{A}$  be a Banach algebra and take  $ m\geq 1$.
Suppose that  $\mathcal{A}$    has pointwise left identity of bound
$ m $ ($i.e.$ for  every $a\in\mathcal{A}$    there exists
  $u\in\mathcal{A}$   with $ \|u\|\leq m $ such that $ua=a$).
Then, for each $ a_{1},...,a_{n} \in \mathcal{A}$ there exists
$u\in\mathcal{A}$    with  $ \|u\|\leq m $  such that $ ua_{i}=a_{i}
,1 \leq i \leq n $.
\begin{proof}
Take $ a_{1},...,a_{n} \in \mathcal{A}$ . Successively choose $
u_{1},...,u_{n} \in \mathcal{A}$ with $\|u_{i}\|\leq m $ , $1\leq
i\leq n $, and $ ( e-u_{i}) ... (e - u_{1})a_{i}=0$  $ (1\leq i \leq
n )$. Here, notice that we use $e$ as a symbol. For instance, by
$(e-u)a$ we mean $a-au$. Define $u\in\mathcal{A}$ by $e-u =(e-u_{n}
)...(e-u_{1})$. Then for each $1\leq i \leq n$ we have
$$a_{i} -ua_{i} = (e-u_{n} )...(e-u_{i+1})((e-u_{i}
)...(e-u_{1})a_{i})=0 \ ,$$ as  required.
\end{proof}
\end{lemma}
\begin{theorem}\label{2.4}
Let   $\mathcal{A}$  be a dual Banach algebra. Suppose that for each
$ a_{1},...,a_{n} \in \mathcal{A}$, there exists
 $u\in\mathcal{A}$      with                  $ \|u\|\leq m $  and    $ ua_{i}=a_{i},$ $1 \leq i \leq n $. Then $\mathcal{A}$      has left identity.
\begin{proof}
Let $\mathcal{F}$   be  the family of all non-empty, finite subset
of  $\mathcal{A}$, ordered by inclusion. Therefore $\mathcal{F}$ is
a directed set. For each $F\in\mathcal{F}$ choose
$e_{F}\in\mathcal{A}$  with
 $\|  e_{F}\|\leq m $ and  $ a= e_{F}a$, $ a\in F$.
 Since  $ (e_{F})_{F\in \mathcal{F}}$ is a  bounded net in a dual Banach algebra, there is
 $e\in\mathcal{A}$  such that $e=\lim e_{F}$ .
  Now, it is clear that $e $ is a left identity for $\mathcal{A}$.
\end{proof}
\end{theorem}

\begin{theorem}\label{2.5}
Let  $\mathcal{A}=(\mathcal{A}_*)^*$  be a commutative pointwise
Connes amenable dual Banach algebra . Then $\mathcal{A}$   has an
identity.
\begin{proof}
By  \cite[Proposition 4.2]{7}, for  every $a\in\mathcal{A}$ there is
a
 bounded net $(u_{i})\subseteq\mathcal{A}$        such that $ au_{i}\xrightarrow{w^{*}} a$    in  $\mathcal{A}$.
   Because   $\mathcal{A}$  is a dual  space, passing to a subnet, we may suppose that there is  $u\in\mathcal{A}$
    such that $u_{i}\longrightarrow u$.
Hence, for every $a\in\mathcal{A}$, there is $u\in\mathcal{A}$ such
that $au=u$. For each $ n=1,2,....$, we define $${A_{n}} =\{a \in
\mathcal{A} : \ ua=a \ \text{for some} \ u \in \mathcal{A} \
\text{with} \  \|u\| \leq n
 \} \ .$$

An argument similar to \cite[Theorem 9.7]{4}, shows  that
$\mathcal{A} $   has pointwise identity of bound $m$
 for  some $m\geq1 $.
Now, Lemma \ref{2.3} and Theorem \ref{2.4} yield that $\mathcal{A} $
possesses an identity.
\end{proof}
\end{theorem}

\section{Examples}
We recall that a \textit{weight} on a discrete semigroup $S$ is a
function $ \omega : S \longrightarrow [1,\infty)$ such that
$\omega(g h) \leq \omega(g) \omega(h)$, for all $g , h \in S$. Then
$\ell^1(S,\omega)$ is the Banach space of all complex functions $ f
= (a_g)_{g \in S}$ on $S$ with the norm $||(a_g)_g || = \sum_{g \in
S} |a_g| \omega(g) < \infty $.

 It is known that $\ell^1(\omega) := \ell^1(\mathbb{N},\omega)$ is a
dual Banach algebra under pointwise multiplication with the predual
$ c_0(\omega^{-1})$ and without identity element, see for instance
\cite{1}.

\begin{theorem}\label{3.1}
Let $\omega$ be a weight on $\mathbb{N}$. Then:

$(i)$ $\ell^1(\omega)$ is not Connes amenable;

$(ii)$ $\ell^1(\omega)$ is not pointwise Connes amenable;

$(iii)$ $\ell^1(\omega)$ is not $w^*$-approximately Connes amenable;

$(iv)$ $\ell^1(\omega)$ is pointwise $w^*$-approximately Connes
amenable.
\begin{proof}
 Since $\ell^1(\omega)$ does not have an identity, the clause $(i)$ is
immediate by \cite[ Proposition 4.1]{9}. The commutativity of
$\ell^1(\omega)$ together with Theorem \ref{2.5} imply $(ii)$. The
clause $(iii)$ is exactly \cite[Theorem 3.3]{7}. Finally, it was
shown in \cite[Corollary 1.8.5]{2}
 that $\ell^1(\omega)$ is pointwise-approximately
amenable. Therefore, automatically it is pointwise
$w^*$-approximately Connes amenable.
\end{proof}
\end{theorem}

Next, we consider the semigroup $\mathbb{N}_{\wedge}$
 which is $\mathbb{N }$  with  the  semigroup  operation
$m \wedge n := min\{m, n\}$,  $(m, n \in \mathbb{N})  $. It is clear
that any function $ \omega : \mathbb{N} \longrightarrow [1,\infty)$
is a weight on the semigroup $\mathbb{N}_{\wedge}$. Then $
\ell^{1}(\mathbb{N}_{\wedge},\omega)$ is a commutative Banach
algebra with the \textit{convolution} product $ \delta_m \star
\delta_n = \delta_{m \wedge n}$, where $\delta_n$ stands for the
characteristic function of $\{n \}$ for $n \in \mathbb{N}$. It is
well known that $\ell^{1}(\mathbb{N}_{\wedge},\omega)$ is a dual
Banach algebra, whenever $\lim_n \omega(n)=\infty$ and with the
predual $ c_0(\omega^{-1})$ \cite[Proposition 3.1.1]{2}.

Suppose that $\lim_n \omega(n)=\infty$. Then, because of the lack of
identity \cite[Propositions 3.3.1 and 3.3.2]{2},
$\ell^{1}(\mathbb{N}_{\wedge},\omega)$ is not Connes amenable. The
same reason, using Theorem \ref{2.5}, implies that
$\ell^{1}(\mathbb{N}_{\wedge},\omega)$ is not pointwise Connes
amenable. It was shown in \cite[Theorem 3.7.1]{2} that
$\ell^{1}(\mathbb{N}_{\wedge},\omega)$ is pointwise approximately
amenable, and therefore it is pointwise $w^*$-approximately Connes
amenable as well. We summarize these facts as follows.
\begin{theorem}\label{3.2}
Let $\omega$ be a function on $\mathbb{N}$ such that $\lim_n
\omega(n)=\infty$. Then:

$(i)$ $\ell^{1}(\mathbb{N}_{\wedge},\omega)$ is not Connes amenable;

$(ii)$ $\ell^{1}(\mathbb{N}_{\wedge},\omega)$ is not pointwise
Connes amenable;

$(iii)$ $\ell^{1}(\mathbb{N}_{\wedge},\omega)$ is pointwise
$w^*$-approximately Connes amenable.

\end{theorem}

 At the end, it should be remarked that $w^*$-approximate Connes amenability
of $\ell^{1}(\mathbb{N}_{\wedge},\omega)$ is an open question for
the authors.

\section{Relations with diagonals }

 Let $\mathcal{A}=(\mathcal{A}_*)^*$ be a dual Banach algebra and let $E$ be a Banach $\mathcal{A}$-bimodule.
  We write $\sigma wc(E)$ for the set of all elements $x\in E$ such that the map
\[\mathcal{A}\longrightarrow E,\quad a\longmapsto \left\{\begin{array}{c}a\cdot x\\ x\cdot a\end{array},\right.\]
are $w^*$-weak continuous. It was shown that
$\pi^*(\mathcal{A}_*)\subseteq \sigma
wc(\mathcal{A}\hat{\otimes}\mathcal{A})^*$  \cite[Corollary
4.6]{11}. Taking adjoints, we can extend $\pi$ to
 an $\mathcal{A}$-bimodule homomorphism $\pi_{\sigma wc}$ from $\sigma wc((\mathcal{A}\hat{\otimes}\mathcal{A})^*)^*$ to
 $\mathcal{A}$.

 Let $E$
be a Banach space. We then have the canonical map $ \imath_E : E
\longrightarrow E^{**}$ defined by $ \langle \mu , \imath_E(x)
\rangle = \langle x , \mu \rangle$ for $\mu \in E^*$, $x \in E$. For
Banach spaces $E$ and $F$, we write $ \mathcal{L}(E,F)$ for the
Banach space of bounded linear maps between $E$ and $F$.  It is
standard that $ (E \hat{\otimes} F)^* = \mathcal{L}(F,E^*)$. For a
Banach algebra $\mathcal{A}$, then we obtain a bimodule structure on
$\mathcal{L}(\mathcal{A},\mathcal{A}^*)= (\mathcal{A} \hat{\otimes}
\mathcal{A})^*$ through $ ( a \ . \ T) (b) = T(ba) \ , \ \ \ (T \ .
\ a)(b) = T(b) \ . \ a $, for $a,b \in \mathcal{A}$, and for $ T \in
\mathcal{L}(\mathcal{A},\mathcal{A}^*)$.

Throughout, we use the term \textit{unital } for a semigroup (or an
algebra) $X$ with an identity element $e_X$.

The following characterizations  will be needed in the sequel.
\begin{theorem}\label{4.1} Let $\mathcal{A}$ be a unital Banach algebra. Then $\mathcal{A}$ is pointwise amenable if and only if for
each $a \in\mathcal{A}$ there exists $ M \in (\mathcal{A}
\hat{\otimes} \mathcal{A})^{**}$ such that $ a \ . \ M = M \ . \ a$,
and $  \pi^{**} (M) = e_{\mathcal{A}}.$
\begin{proof}
The proof is a small variation of the standard argument in
\cite[Theorem 43.9]{0}.
\end{proof}
\end{theorem}
\begin{theorem}\label{4.2} Let $\mathcal{A}$ be a unital dual Banach algebra.
Then:

$(i)$ $\mathcal{A}$ is pointwise Connes amenable if and only if for
each $a \in\mathcal{A}$ there exists $M \in \sigma wc((\mathcal{A}
\hat{\otimes} \mathcal{A})^*)^*$ such that $ a \ . \ M = M \ . \ a$,
and $  \pi_{\sigma wc} (M) = e_{\mathcal{A}}.$

$(ii)$ $\mathcal{A}$ is pointwise Connes amenable if and only if for
each $a \in\mathcal{A}$ there exists $ M \in (\mathcal{A}
\hat{\otimes} \mathcal{A})^{**}$ such that $ \langle T , a \ . \ M -
M \ . \ a \rangle = 0$ for each $T \in \sigma wc ( \mathcal{L}(
\mathcal{A} , \mathcal{A}^*))$, and $\imath_{\mathcal{A}_*}^*
\pi^{**} (M) =e_{\mathcal{A}}.$
\begin{proof}
The clause $(i)$ is analogous to \cite[Theorem 4.8]{11}. Because
$\sigma wc((\mathcal{A} \hat{\otimes} \mathcal{A})^*)^*$ is a
quotient of $(\mathcal{A} \hat{\otimes} \mathcal{A})^{**}$, the
clause $(ii)$ is just a re-statement of $(i)$.
\end{proof}
\end{theorem}

For a discrete semigroup $S$, we recall that $ \ell^1(S)
\hat{\otimes} \ell^1(S) = \ell^1(S \times S) $, where $ \delta_g
\otimes \delta_h$ is identified with $\delta_{(g,h)}$ for $ g , h
\in S$. Thus we have $ \mathcal{L}(\ell^1(S) , \ell^\infty(S) ) = (
\ell^1(S) \hat{\otimes} \ell^1(S))^* = \ell^1(S \times S)^* =
\ell^\infty(S \times S) $, where $ T \in \mathcal{L}(\ell^1(S) ,
\ell^\infty(S) )$ is identified with $ (T_{(g,h)})_{(g,h) \in S
\times S} \in \ell^\infty(S \times S)$, while $ T_{(g,h)} := \langle
\delta_h , T(\delta_g) \rangle$. Let $\omega$
 be a weight on $S$. If $S$ is unital then, without
loss of generality, we put $ \omega(e_S) = 1$. The Banach space
$\ell^1(S,\omega) $ with the convolution product is a Banach
algebra, called a \textit{Beurling} algebra. We consider
$\ell^1(S,\omega)$ as the Banach space $\ell^1(S)$ with the product
$ \delta_g \star_\omega \delta_h := \delta_{gh} \Omega(g,h)$, where
$\Omega(g,h) := \frac{\omega(gh)}{\omega(g) \omega(h)}$, $(g,h \in
S)$, and extend $\star_\omega$ to $\ell^1(S)$ by linearity and
continuity. A semigroup $S$ is \textit{weakly cancellative} if, for
each $s \in S$, the maps $L_s$ and $R_s$, defined by $ L_s(t) = st$
and $R_s(t)= ts$, are finite-to-one. In this case,
$\ell^1(S,\omega)$ is a dual Banach algebra with predual $c_0(S)$
\cite[Proposition 5.1]{12}

\begin{theorem}\label{4.3} Let $S$ be a discrete unital semigroup, let $\omega$ be a weight on $S$ and let $ \mathcal{A}:=
\ell^1(S,\omega)$. Consider the following statements:

$(1)$ $\mathcal{A}$ is pointwise amenable.

$(2)$ For each $k \in S$ there exists $ M \in (\mathcal{A}
\hat{\otimes} \mathcal{A})^{**}= \ell^\infty(S \times S)^*$ such
that:

$(i)$ $\langle (f(hk,g) \Omega(h,k) - f(h,kg) \Omega(k,g))_{(g,h)
\in S \times S} , M \rangle = 0 $ for each bounded function $f : S
\times S \longrightarrow \mathbb{C} $;

$(ii)$ $ \langle (f_{gh} \Omega(g,h))_{(g,h) \in S \times S} , M
\rangle = f_{e_S}$ for each bounded family $(f_g)_{g \in S}$.

$(3)$ $\mathcal{A}$ is pointwise amenable at $\delta_k$ for each $k
\in S$.

Then we have $(1) \Longrightarrow (2) \Longrightarrow (3)$.
\begin{proof}
$(1) \Longrightarrow (2)$ Suppose that $\mathcal{A}$ is pointwise
amenable and that $k \in S$. Take $M \in \ell^\infty(S \times S)^*$
as in Theorem \ref{4.1}. For each bounded family $(f_g)_{g \in S}$,
we have
$$ \pi^*(f) = ( \langle \delta_{gh}, f \rangle \Omega(g,h) )_{(g,h)
\in S \times S} \in \ell^\infty(S \times S) \ .$$ Therefore
$$  \langle (f_{gh} \Omega(g,h))_{(g,h) \in S \times S}, M \rangle
=  \langle  f, \pi^{**}(M) \rangle  =  \langle  f, e_{\mathcal{A}}
\rangle =  f_{e_S} .$$ Next, for every $T \in \mathcal{L}(
\mathcal{A}, \mathcal{A}^* ) =\ell^\infty(S \times S)$, we see that
$$ \langle \delta_g \otimes \delta_h , \delta_k \ . \ T -  T \ . \
\delta_k \rangle = \langle \delta_g , T (\delta_{hk}) \rangle
\Omega(h,k) - \langle \delta_{kg} , T (\delta_h) \rangle \Omega(k,g)
\ .$$ Let $f : S \times S \longrightarrow \mathbb{C} $ be a bounded
function and take $T \in \mathcal{L}( \mathcal{A}, \mathcal{A}^* )
=\ell^\infty(S \times S)$ defined by $ \langle \delta_h ,
T(\delta_g) \rangle = f(g,h)$. Hence
\begin{align*}
 \langle (f(hk,g) \Omega(h,k) - f(h,kg) \Omega(k,g))_{(g,h)
\in S \times S} , M \rangle  &=  \langle ( \langle \delta_g \otimes
\delta_h , \delta_k \ . \ T -  T \ . \ \delta_k \rangle)_{(g,h) \in
S \times S} , M \rangle    \\&= \langle \delta_k \ . \ T -  T \ . \
\delta_k , M  \rangle
\\&= \langle T , \delta_k \ . \ M -  M \ . \ \delta_k  \rangle = 0,
\end{align*}
as required. Likewise for the implication $(2) \Longrightarrow (3)$.
\end{proof}
\end{theorem}
 The following is \cite[Proposition 5.5]{12}, in which $\mathcal{W}(\mathcal{A}, \mathcal{A}^*
 )$ stands for the collection of weakly compact operators in $\mathcal{L}( \mathcal{A}, \mathcal{A}^*
 )$, and the set of weakly almost periodic elements in $\mathcal{W}(\mathcal{A}, \mathcal{A}^*
 )$ is denoted by $WAP
(\mathcal{W}(\mathcal{A}, \mathcal{A}^* ))$.

\begin{theorem}\label{4.4} Let $S$ be a weakly cancellative semigroup, let
$\omega$ be a weight on $S$, and let $\mathcal{A}:=
\ell^1(S,\omega)$. Let $T \in \mathcal{L}( \mathcal{A},
\mathcal{A}^* )$ be such that $ T(\mathcal{A}) \subseteq
\imath_{c_0(S)} ( c_0(S))$ and $T^*(\imath_{\mathcal{A}}
(\mathcal{A})) \subseteq \imath_{c_0(S)} ( c_0(S))$. Then $ T \in
\mathcal{W}(\mathcal{A}, \mathcal{A}^* )$, and $ T \in WAP
(\mathcal{W}(\mathcal{A}, \mathcal{A}^* ))$ if and only if, for each
sequence $(k_n)$ of distinct elements of $S$, and each sequence
$(g_m,h_m)$ of distinct elements of $S \times S$ such that the
repeated limits
$$ \lim_n \lim_m \langle \delta_{k_n g_m} , T(\delta_{h_m})  \rangle
\ \ , \ \ \lim_n \lim_m \Omega(k_n, g_m) $$
$$ \lim_n \lim_m \langle \delta_{ g_m} ,
T(\delta_{h_m k_n})  \rangle \ \ , \ \ \lim_n \lim_m \Omega(h_m ,
k_n) $$ all exist, we have at least one repeated limit in each row
is zero.
\end{theorem}

\begin{theorem}\label{4.5} Let $S$ be a discrete, weakly cancellative semigroup, let $\omega$ be a weight on $S$ and let $ \mathcal{A}:=
\ell^1(S,\omega)$ be unital.  Consider the following statements:

$(1)$ $\mathcal{A}$ is pointwise Connes amenable.

$(2)$ For each $k \in S$ there exists $ M \in (\mathcal{A}
\hat{\otimes} \mathcal{A})^{**}= \ell^\infty(S \times S)^*$ such
that:

$(i)$ $\langle (f(hk,g) \Omega(h,k) - f(h,kg) \Omega(k,g))_{(g,h)
\in S \times S} , M \rangle = 0 $ for each bounded function $f : S
\times S \longrightarrow \mathbb{C} $ which is such that the map $T
\in \mathcal{L}(\mathcal{A}, \mathcal{A}^*)$, defined by $ \langle
\delta_h , T(\delta_g) \rangle = f(g,h)$, for $g,h \in S$, satisfies
the conclusions of Theorem \ref{4.4};

$(ii)$ $ \langle (f_{gh} \Omega(g,h))_{(g,h) \in S \times S} , M
\rangle = \langle f, e_{\mathcal{A}} \rangle $ for each family
$(f_g)_{g \in S} \in c_0(S)$.

$(3)$ $\mathcal{A}$ is pointwise Connes amenable at $\delta_k$ for
each $k \in S$.

Then we have $(1) \Longrightarrow (2) \Longrightarrow (3)$.
\begin{proof}
Using Theorem \ref{4.2} $(ii)$ in place of Theorem
 \ref{4.1}, this follows as Theorem \ref{4.3}.
\end{proof}
\end{theorem}

Let $G$ be a discrete group and let $h \in G$. Define $ J_h :
\ell^\infty(G) \longrightarrow \ell^\infty(G)$ by
$$ J_h(f) := ( f_{hg} \Omega(h,g) \omega(h) \Omega(g^{-1}, h^{-1})
\omega(h^{-1}))_{g \in G} \ \ \ ( f=(f_g)_g \in \ell^\infty(G)) \
.$$ It is clear that $ ||J_h(f)|| \leq ||f||
\omega(h)\omega(h^{-1})$, so that $J_h$ is bounded.

The following is the pointwise variant of \cite[Definition
5.10]{12}.
\begin{definition}\label{4.6} Let $\omega$ be a weight on a discrete group $G$,
and let $h_0 \in G$. Then $G$ is \textit{pointwise}
$\omega$-\textit{amenable at} $h_0$ if there exists $N \in
\ell^\infty(G)^*$ such that $ \langle (\Omega(g, g^{-1})_{g \in G},
N \rangle = 1 $ and $ J_{h_0}^* (N) = N.$ We say that $G$ is
\textit{pointwise} $\omega$-\textit{amenable} if it is pointwise
$\omega$-amenable at each $h \in G$.
\end{definition}

\begin{theorem}\label{4.7} Let $G$ be a discrete group, let $\omega$ be a weight on $G$  and let $\mathcal{A} =
\ell^1(G,\omega)$. Consider the following statements:

$(1)$ $\mathcal{A}$ is pointwise amenable.

$(2)$ $ \mathcal{A}$ is pointwise Connes-amenable.

$(3)$  $G$ is pointwise $\omega$-amenable.

$(4)$ $\mathcal{A}$ is pointwise amenable at $\delta_k$ for each $k
\in G$.

Then we have $(1) \Longrightarrow (2) \Longrightarrow (3)
\Longrightarrow (4)$.
\begin{proof}
The implication  $(1) \Longrightarrow (2) $ is trivial, and $(2)
\Longrightarrow (3)$ is a more or less verbatim of the proof of
\cite[Theorem 5.11 $(1) \Longrightarrow (3)$]{12} .

$(3) \Longrightarrow (4)$: Let $k \in G$, and let $N  \in
\ell^\infty(G)^*$ be given as in $(3)$. Define $ \psi :
\ell^\infty(G \times G) \longrightarrow \ell^\infty(G)$ by $ \langle
\delta_g , \psi(F) \rangle := F(g,g^{-1})$, for each $ F \in
\ell^\infty(G \times G)$ and $ g \in G$. Put $ M := \psi^*(N)$. Then
it suffices to show that
 $M$ has desired properties in Theorem \ref{4.3}
$(2)$. First, for a bounded family $(f_g)_{g\in G}$, we see that $$
\langle (f_{gh} \Omega(g,h))_{(g,h)} , M \rangle  = \langle (f_{e_G}
\Omega(g,g^{-1}))_g , N  \rangle=  \langle ( \Omega(g,g^{-1}))_g , N
\rangle= 1 \ .
$$

Next, for an arbitrary bounded function $f :G \times G
\longrightarrow \mathbb{C}$ , it is clear that $$ \psi((f(hk,g)
\Omega(h,k)- f(h,kg) \Omega(k,g))_{(g,h)}) = ( f(g^{-1}k,g)
\Omega(g^{-1},k) - f(g^{-1},kg) \Omega(k,g))_g \ .$$ Define $ F: G
\times G \longrightarrow \mathbb{C}$, by $ F(g,h) :=f(hk,g)
\Omega(h,k)$, for each $g,h \in G$. Hence, it is readily seen that
$F$ is bounded and $ ||F||_\infty \leq ||f||_\infty$. Therefore
\begin{align*}
 \langle ( f(hk,g) \Omega(h,k) - f(h,kg) \Omega(k,g))_{(g,h)} ,
M \rangle  &=  \langle ( f( g^{-1}k, g) \Omega(g^{-1},k) - f(g^{-1}
, kg) \Omega(k,g) )_g , N \rangle  \\&=  \langle \psi(F) - J_k (
\psi(F)) , N \rangle \\&=  \langle \psi(F) , N - J_k^* ( N) \rangle
= 0  ,
\end{align*}
as required.
\end{proof}
\end{theorem}
\begin{remark}\label{4.8} It seems to be a \textit{right} conjecture that pointwise amenability of $\ell^1(G,\omega)$ coincides
 with its pointwise Connes amenability, however we are not able to prove
 (or disprove) it. In fact, we can not establish the implication $(4) \Longrightarrow (1)$ in Theorem \ref{4.7},
 because we can see no reason that pointwise amenability at elements $a$ and $b$ gives any information about $a+b$.
\end{remark}



\end{document}